\documentclass[12pt,a4paper,twoside,final,notitlepage, leqno]{article}
\usepackage[english]{babel}
\usepackage[T1]{fontenc}  
\usepackage{epsfig, graphicx, amssymb}
\usepackage{amsmath,amsthm,epsfig,amsfonts}  
\usepackage{float}
\usepackage{color}
\def\br {\break}

\def \smb {{\scriptstyle \bullet }}
\newcommand{\monitem}{ \smallskip \noindent $\bullet$ \quad  } 
\newcommand{\moneq}{\vspace*{-7pt} \begin{equation} \displaystyle } 
\newcommand{\moneqstar}{\vspace*{-6pt} \begin{equation*} \displaystyle } 
\newcommand{\monendstar}{\vspace*{-6pt} \end{equation*}   }
\newcommand{\monend}{\vspace*{-7pt} \end{equation}   }

\newcommand{\dd}{{\rm d}}
\newcommand{\R}{\mathbb{R}}

\newcommand{\RR}[0]{\mathbb{R}}

\def\section*#1{}

\usepackage{fancyhdr}
\renewcommand{\headrulewidth}{0pt}

\begin{document} 

\fancypagestyle{plain}{ \fancyfoot{} \renewcommand{\footrulewidth}{0pt}}
\fancypagestyle{plain}{ \fancyhead{} \renewcommand{\headrulewidth}{0pt}}

~

  \vskip 2.1 cm

  \centerline {\bf \LARGE  Equivalent partial differential equations}
  
 \bigskip 

\centerline {\bf \LARGE  of a lattice Boltzmann scheme }

 \bigskip  \bigskip \bigskip

\centerline { \large    Fran\c{c}ois Dubois$^{ab}$ }

\smallskip  \bigskip 

\centerline { \it  \small   
$^a$   Dpt. of Mathematics, University Paris-Sud,  B\^at. 425, F-91405  Orsay, France.} 

\centerline { \it  \small   
$^b$    Conservatoire National des Arts et M\'etiers, LMSSC laboratory,  F-75003 Paris, France.}


\bigskip  

\centerline {  05 septembre 2007  
  {\footnote {\rm  \small $\,$ Contribution published in 
 {\it Computers \& Mathematics with Applications}, volume 55, pages~1441-1449,
 2008, doi https://doi.org/10.1016/j.camwa.2007.08.003. Edition  07 June 2018. }} }


 
 {\bf AMS (MSC2010) classification}: 76M28.

\bigskip  
\noindent {\bf \large Abstract}

 \noindent
We show that when we formulate the lattice Boltzmann
equation with a small time step  $ \, \Delta t \, $ and an 
associated space scale $ \, \Delta x , \, $ a 
Taylor expansion joined with the so-called equivalent equation methodology leads  to 
establish macroscopic fluid equations as a formal limit. We recover the Euler equations of
gas dynamics at the first order and the compressible  Navier-Stokes  
equations at the second order.

\bigskip \bigskip    \noindent {\bf \large    1) \quad  Discrete geometry }    

\monitem 
We denote by $d$ the dimension of  space and by 
 $\, \cal{ L} \, $   a regular $d$-dimensional lattice. Such a lattice 
  is composed by a set  $\, \cal{ L}^{\rm 0} \, $ of nodes or vertices and a set 
 $\, \cal{ L}^{\rm 1} \, $ of links or edges between two vertices. 
From a practical point of view, given
 a vertex $\, x ,\,$ there exists a set  $\, V(x) \,$ of neighbouring nodes, including the
 node  $\, x \,$ itself. 
We consider here that the   lattice   $\, \cal{ L} \, $ is parametrized by a space 
step  $\, \Delta x  > 0 . \, $ 
For the  fundamental example called D2Q9 (see {\it e.g.} Lallemand and Luo, 2000), 
the set  $\, V(x) \,$ is given with the help of the family of vectors $\, (e_j)_{0 \leq j
\leq J} \,$ defined by $\, J=8 ,\,$ 

\setbox21=\hbox{$\displaystyle 
  \begin{pmatrix} 0 \\ 0 \end{pmatrix} ,\,
  \begin{pmatrix} 1 \\ 0 \end{pmatrix} ,\,
  \begin{pmatrix} 0 \\ 1 \end{pmatrix} ,\,
  \begin{pmatrix} -1 \\ 0 \end{pmatrix} ,\,
  \begin{pmatrix} 0 \\ -1 \end{pmatrix} ,\,
  \begin{pmatrix} 1 \\ 1 \end{pmatrix} ,\,
  \begin{pmatrix} -1 \\ 1 \end{pmatrix} ,\,
  \begin{pmatrix} -1 \\ -1 \end{pmatrix} ,\,
  \begin{pmatrix} 1 \\ -1 \end{pmatrix}  $}
\setbox31= \hbox{ $\vcenter {\box21} $}
\setbox44=\hbox{\noindent  (1.1) $\displaystyle  \qquad  \big(  e_j  \big) 
\, = \,  \left\{ \box31 \right\} $}  
\smallskip \noindent $ \box44 $

\smallskip \noindent
and  the vicinity

\smallskip \noindent   (1.2)  $\qquad   \displaystyle
V(x) \, = \, \{ \, x \,+\, \Delta x \, e_j, \,\, 0 \leq j \leq J \, \}  \, . \, $ 

\bigskip \noindent $\bullet$ \quad 
In the general case, we still suppose that the equation 
(1.2) holds but we do not make any precise
definition concerning the integer  $\, J \,$  and the nondimensionalized vectors 
 $\, (e_j)_{0 \leq j \leq J}.  \,$ Nevertheless  if $\, x \, $ is a node of the lattice 
($ x \in  \cal{ L}^{\rm 0} ) $, then $ \,  y^j  =  x \,+\, \Delta x \, \,  e_j  \, $ 
is  an other node of the lattice, {\it i.e.} $ \, y^j  \in   \cal{ L}^{\rm 0}. \, $

\bigskip \bigskip    \noindent {\bf \large    2) \quad  Lattice Boltzmann framework }    

\monitem 
We introduce a time step  $\, \Delta t > 0 \, $ and we suppose that the celerity $\,
\lambda \, $ defined according to 

\smallskip \noindent   (2.1)  $ \qquad   \displaystyle
\lambda    \, =   \, {{\Delta x}\over{ \Delta t }} \, $

\smallskip \noindent
remains  fixed. Then we introduce a local velocity $\, v_j \,$ in such a way that 

\smallskip \noindent   (2.2)  $ \qquad   \displaystyle
  \Delta t  \, \,  v_j \, \, = \, \Delta x \,  \, e_j   \,, \quad  0 \leq j \leq J  \,. \, $ 

\smallskip \noindent 
In this $d$-dimensional framework we will denote by  $\, v_j^{\alpha} \, ( 1 \leq \alpha
\leq d)     \,$ the Cartesian components of velocities  $\, v_j .\,$ Recall that if 
 $\, x \,$ is a node of the lattice,  the point   $\,  x \,+\,  \Delta t  \, \,  v_j  \,$ 
is also a node of the lattice:  

\smallskip \noindent   (2.3)  $ \qquad   \displaystyle
 x \in  {\cal{ L}}^{\rm 0} \, \, \Longrightarrow \,  \,   x \,+\,  \Delta t  \, \,  v_j
 \in {\cal{ L}}^{\rm 0} \,, \quad \forall j = 0, \dots J. \, $

\bigskip \noindent $\bullet$ \quad 
According to D'Humi\`eres (1992), the lattice Boltzmann
scheme describes the dynamics of the density $\, f^j(x,\,t) \, $  of particles of velocity 
 $\, v_j \,$ at the node $ \, x \, $ and for  the discrete time  $ \, t . \, $ 
We introduce the  $ \, d+1 \, $ scalar ``conservative variables''  $ \, W(x,\, t) \, $
composed by the density $ \, \rho \, $ and the momentum $ \, q .\, $ Note that it is also
possible to take into account the conservation of the total energy (see D'Humi\`eres's
article  for example). We have 

\noindent   (2.4)  $ \qquad   \displaystyle
\rho(x,\, t)  \, = \, \sum_{j=0}^{J} \, f^j (x,\, t) \,\equiv\, W^0 (x,\, t) \, $ 

\smallskip \noindent   (2.5)  $ \qquad   \displaystyle
q^{\alpha} (x,\, t)  \, = \, \sum_{j=0}^{J} \, v_j^{\alpha}
 \,  f^j (x,\, t) \,\equiv\, W^{\alpha} (x,\, t) \,, \quad  1 \leq \alpha \leq d  \,, \, $ 

 \noindent 
 and
 
\smallskip \noindent   (2.6)  $ \qquad   \displaystyle
 W(x,\, t)   \, = \, \big( \, \rho(x,\, t) ,\, q^{1} (x,\, t) ,\, \cdots \, ,\, 
  q^{d} (x,\, t) \, \big) \, . \, $ 

\smallskip \noindent 
When a state $ \,  W  \, $ is given in space $\, \R^{d+1} ,\, $ a Gaussian (or any other
choice)  equilibrium distribution of particles is defined according to 

\smallskip \noindent   (2.7)  $ \qquad   \displaystyle
f^j_{\rm eq}    \, = \, G^j (W)  \,, \quad  0 \leq j \leq J  \, \, $ 

\smallskip \noindent  
in such a way that 

\setbox44=\hbox{\noindent  (2.8) $\displaystyle \qquad  \sum_{j=0}^{J} \, G^j (W)  \, \equiv \, W^0   \,, \quad
\sum_{j=0}^{J} \, v_j^{\alpha}  \,  G^j (W)
\, \equiv\, W^{\alpha} \,, \quad  1 \leq \alpha \leq d  \,. $} 
\smallskip \noindent $ \box44 $

\fancyhead[EC]{\sc{Fran\c cois  Dubois }} 
\fancyhead[OC]{\sc{Equivalent partial differential equations of a lattice Boltzmann scheme}} 
\fancyfoot[C]{\oldstylenums{\thepage}}

\bigskip \noindent $\bullet$ \quad 
Following D'Humi\`eres (1992), we introduce the ``moment vector'' $ \, m \, $ according to 

 \noindent   (2.9)  $ \qquad   \displaystyle
 m^{k} \, =  \,   \sum_{j=0}^{J}  \, M^k_j \,  \,  f^j  \,, \quad  0 \leq k  \leq J  \,. \, $ 

\smallskip \noindent 
For $ \, 0 \leq i \leq d ,\, $ the moments $\, m^i \, $ are {\bf identical} to the
conservative variables: 

\smallskip \noindent   (2.10)  $ \qquad   \displaystyle
m^0  \, \equiv \,  \rho \,, \quad m^{\alpha}  \, \equiv \,  q^{\alpha} 
 \,, \quad  1 \leq \alpha \leq d  \,. \, $ 

\smallskip \noindent 
In other words, the matrix $ \, M \, $ satisfies 

\smallskip \noindent   (2.11)  $ \qquad   \displaystyle
M^0_j  \, \equiv \, 1  \,, \quad M^{\alpha}_j  \, \equiv \,  v^{\alpha}_j  
\,, \quad  0 \leq j \leq J  \,,  \quad  1 \leq \alpha \leq d  \,. \, $ 

\smallskip \noindent
We assume that vectors  $\, (e_j)_{0 \leq j \leq J} \,$  are chosen such that the 
$ \, (d+1) \times (J+1) \,$ matrix 
$ \, (M_{k \, j})_{0 \leq k \leq d ,\,  0 \leq j \leq J}  \,$
is of full rank. With this hypothesis, the conservative moments $ \, W \, $ introduced in relation
(2.6) are independent variables. 

\bigskip \noindent $\bullet$ \quad 
When a particle distribution $ \, f \, $ is given, the moments are evaluated according to
(2.9). The matrix  $ \, M \, $ is supposed to be invertible and the inverse relation takes
the form: 

\smallskip \noindent   (2.12)  $ \qquad   \displaystyle
  f^j  \, =  \,   \sum_{k=0}^{J}  \, (M^{-1})^j_k   \,  \,  m^{k}  
\,, \quad  0 \leq j  \leq J  \,. \, $ 

\smallskip \noindent 
When $ \,   f^j_{\rm eq}  \, $ is determined according to the relation (2.7), the
associated equilibrium  moments $ \,  m^{k}_{\rm eq}  \, $ are given simply according 
to (2.9), {\it i.e.} in this case 

\smallskip \noindent   (2.13)  $ \qquad   \displaystyle
 m^{k}_{\rm eq}  \, =  \,   \sum_{j=0}^{J}  \, M^k_j \,  \,  f^j_{\rm eq} 
  \,, \quad  0 \leq k  \leq J  \,. \, $

\smallskip \noindent  
We remark also that by construction (relation (2.8)), we have 

\smallskip \noindent   (2.14)  $ \qquad   \displaystyle
 m^{i}_{\rm eq}  \, = \,  m^i  \, =  \,  W^i    \,, \quad  0 \leq i \leq d  \,. \, $

\bigskip \bigskip    \noindent {\bf \large    3) \quad Collision step }    

\monitem 
The collision step is local in space and is naturally defined in the space of moments. If 
$ \,  m^{k} (x,\, t) \, $ denotes the value of the $ \, k^{\rm th} \, $ component of the
moment vector $\,  m\, $ at position  $\,  x \, $ and time  $\,  t , \, $ the same component 
$ \,  m^{k}_{*} (x,\, t) \, $ of the moment {\bf after} the collision is 
trivial by construction for the conservative variables: 

\smallskip \noindent   (3.1)  $ \qquad   \displaystyle
 m^{i}_{*} (x,\, t) \,= \,  m^{i} (x,\, t)   \,, \quad  0 \leq i \leq d  \,. \, $ 

\smallskip \noindent 
For the non-conservative components of the moment vector, we {\bf fix} the ratio 
$\, s_k \, $  ($ k \geq d+1 $) between the time step $\, \Delta t \, $ and the relaxation
time $\, \tau_k \, $ of an underlying process: 

\smallskip \noindent   (3.2)  $ \qquad   \displaystyle
s_k \, = \, {{ \Delta t }\over{  \tau_k}}  \,, \quad  d+1 \leq k \leq J  \,. \, $

\bigskip \noindent $\bullet$ \quad  
Then $ \,  m^{k}_{*} (x,\, t) \,$ after the collision is defined according to 

\smallskip \noindent   (3.3)  $ \qquad   \displaystyle
  m^{k}_{*} (x,\, t) \,=\, (1 - s_k) \,   m^{k} (x,\, t)  \, + \,  s_k  \,   m^{k}_{\rm eq} 
 \,, \quad  d+1 \leq k \leq J  \,. \, $ 

\bigskip \noindent  {\bf Proposition 1}. \quad Explicit Euler scheme. 

 \noindent
The numerical scheme (3.3) is exactly the explicit Euler scheme relative to the 
continuous in time relaxation equation 

\smallskip \noindent   (3.4)  $ \qquad   \displaystyle
{{\rm d}\over{{\rm d}t}} \big( m^k - m^k_{\rm eq} \big) 
\, + \, {{1 }\over{  \tau_k}}  \big( m^k - m^k_{\rm eq} \big) \, = \, 0 
 \,, \quad  d+1 \leq k \leq J  \,. \, $ 

 \smallskip
 \newpage \noindent  {\bf Proof of Proposition 1}.

 \noindent
Following {\it e.g.} Strang (1986), we know that the explicit Euler scheme for the
evolution (3.4) takes the form 

\smallskip \noindent   (3.5)  $ \quad   \displaystyle
 {{1 }\over{   \Delta t }} \Big[ \big( m^k - m^k_{\rm eq} \big) (t+\Delta t) - 
\big( m^k - m^k_{\rm eq} \big) (t) \Big] 
\, + \, {{1 }\over{  \tau_k}}  \big( m^k - m^k_{\rm eq} \big)(t) \, = \, 0  \, . \, $ 

\smallskip \noindent
We have by construction the relation (3.1), that is  $\, m^{i} (t+  \Delta t) = 
 m^{i} (t) \, $ for $\,  0 \leq i \leq d \,$ with these notations. Then 
$\, W(t+  \Delta t) = W(t) \,$ and, due to the relation (2.7), 
$\, f^j_{\rm eq} (t+  \Delta t) = f^j_{\rm eq} (t) \,$ after the collision step for all the
components $\, j \, $  of the particle distribution. 
Due to (2.13), we deduce that 
$\,  m^k_{\rm eq}  (t+\Delta t) = m^k_{\rm eq}  (t) \,$ for all  $\,  k \leq J .\,$ 
Thus the expression (3.5) takes the simpler form 

\smallskip \noindent   (3.6)  $ \quad   \displaystyle
 {{1 }\over{   \Delta t }} \big[  m^k(t+\Delta t) -  m^k(t) \big] 
\, + \, {{1 }\over{  \tau_k}}  \big( m^k - m^k_{\rm eq} \big)(t) \, = \, 0  \, , \, $ 

\smallskip \noindent  
which is exactly (3.3), except the change of notations: $\, m^k(t+\Delta t) \,$ 
is replaced by   $\,    m^{k}_{*}  .\,$ $\hfill \square$

\bigskip \noindent $\bullet$ \quad 
We remark also that the classical stability condition for the explicit Euler scheme (see
again  {\it e.g.} the book of Strang) takes the form 

\smallskip \noindent   (3.7)  $ \qquad   \displaystyle
0 \leq \Delta t \leq 2 \, \tau_k \,.\, $ 

\smallskip \noindent 
We will suppose in the following that 

\smallskip \noindent   (3.8)  $ \qquad   \displaystyle
0 \, < \, s_k \, \leq \, 2 \,,   \quad  d+1 \leq k \leq J  \,. \, $ 

\smallskip \noindent  
to put in  
evidence that the moments $ \, m^k \,$ are {\bf not} conserved for index $k$ greater than 
$ d+1 .$ 
We remark also that for the physically relevant Boltzmann equation, the relaxation times
$\, \tau_k \,$  have
a physical sense. With the lattice Boltzmann scheme itself, these physical constants are
no longer  correctly approximated     whereas the {\bf  ratios}  
$ \, s_k =  {{\Delta t}\over{\tau_k}} \, $ are
supposed to be fixed in all what follows. Despite the usual ``LBE'' denomination, a lattice
Boltzmann scheme is not a numerical method to approach the Boltzmann equation ! 

\bigskip \noindent $\bullet$ \quad 
The particle distribution $\, f^j _* \, $ after the collision step follows the relation
(2.12). We have precisely after the collision step 

\smallskip \noindent   (3.9)  $ \qquad   \displaystyle
  f^j_*  \, =  \,   \sum_{k=0}^{J}  \, (M^{-1})^j_k   \,  \,  m^{k}_*  
\,, \quad  0 \leq j  \leq J  \,. \, $

\bigskip \bigskip    \noindent {\bf \large    4) \quad Advection step }    

\monitem 
The avection step of the lattice Boltzmann scheme claims that after the collision step,
the particles having  velocity  $ \, v_j \,$ at  position $ \, x \, $ go in one
time step $\, \Delta t \, $ to the $\, j^{\rm th} \, $ neighbouring vertex. Thus the
particle density $\, f^j (x+ v_j \,\Delta t ,\, t+ \Delta t) \, $ at the new time step 
in the neighbouring vertex is equal to the previous particle density 
 $\, f^j_*  (x ,\, t ) \, $ at the position  $ \, x \, $  after the collision: 

\smallskip \noindent   (4.1)  $ \qquad   \displaystyle
 f^j (x+ v_j \,\Delta t ,\, t+ \Delta t) \,= \, f^j_*  (x ,\, t ) \, . \, $ 

\smallskip \noindent 
We re-write this relation in term of the ``arrival'' node $\, x+ v_j \,\Delta t .\, $ 
We set $ \, \widetilde{x} =  x+ v_j \,\Delta t , \, $ then we have 
 $ \,x = \widetilde{x} -  v_j \,\Delta t  \, $ and going back to the notation 
 $\, x , \, $ we write the relation (4.1) in the equivalent manner 

\smallskip \noindent   (4.2)  $ \qquad   \displaystyle
 f^j (x ,\, t+ \Delta t) \,= \, f^j_*  (x -  v_j \,\Delta t ,\, t ) 
\,, \quad  0 \leq j  \leq J  \,, \quad x \in {\cal{L}}^0\,. \, $

\bigskip \noindent  {\bf Proposition 2}. \quad Upwind scheme for the advection equation. 

 \noindent 
The scheme (4.2) for the advection step of the lattice Boltzmann method is nothing else
that the explicit upwind scheme for the advection equation 

\smallskip \noindent   (4.3)  $ \qquad   \displaystyle
{{ \partial f^j }\over{ \partial t }} \, + \,  v_j \, \smb \, \nabla  f^j \, = \, 0 
\,, \quad  0 \leq j  \leq J  \,, \, $ 

\smallskip \noindent 
with a so-called Courant-Friedrichs-Lewy number $\, \sigma_j \,$ in the $\, j^{\rm th} \, $
direction of the lattice defined by 

\smallskip \noindent   (4.4)  $ \qquad   \displaystyle
\sigma_j \,\, \equiv \,\, \mid  v_j \mid  {{ \Delta t }\over{ \Delta x \, \mid  e_j \mid }} \,$

\smallskip \noindent  
equal, due to the definition (2.2),  to  unity: $ \, \sigma_j  = 1 .\, $

\smallskip \noindent  {\bf Proof of Proposition 2}.

 \noindent
When the Courant-Friedrichs-Lewy number  $\, \sigma_j \,$ is equal to unity, it is
classical (see {\it e.g.} Strang, 1986) that the upwind scheme is exact for the advection
equation.   $\hfill \square$

\bigskip \bigskip    \noindent {\bf \large    5) \quad Equivalent equation at zero order  }    

\monitem
 The lattice Boltzmann scheme is defined by the relations (2.4) to (2.9), (3.3) and
 (4.2). It is parametrized by the lattice step $\, \Delta x ,\, $ the matrix $\, M \, $
 linking the particle distribution $\, f \, $ and the moment vector  $\, m , \, $ the choice of
 the conservative moments, the nonlinear equilibrium function $\, G(\smb) ,\,$ the time step 
 $\, \Delta t \, $ and the ratios  $\, s_k \, $ between the time step and the collision
 time constants for nonequilibrium moments. In what follows, we fix the geometrical and
 topological structure of the lattice $\, \cal{L} ,\, $ we fix the matrix  $\, M \, $
and the equilibrium function  $\, G(\smb) ,\,$ we fix also the ratio  $\, \lambda \, $
defined in (2.1) and last but not least, we suppose that the parameters   $\, s_k \, $
for $\, k \geq d+1 \,$ have a fixed value. Then the whole lattice Boltzmann scheme depends
 on the single parameter $ \, \Delta t. \, $

\bigskip \noindent $\bullet$ \quad
We explore now formally what are the partial 
differential equations associated with the Boltzmann   numerical scheme, following the so-called
 ``equivalent equation method'' introduced and developed by 
 Lerat-Peyret (1974) and Warming-Hyett (1974).
 This approach is based on the assumption, that a sufficiently
smooth function exists which satisfies the difference equation at the grid
points. This assumption gives formal responses to put in evidence
partial differential equations that minimimize the   truncation errors  of
 the numerical scheme. Nevertheless, we note here that this method of analysis
 fails to predict initial layers and
boundary effects properly, as discussed by   Griffiths   
and Sanz-Serna (1986) or  Chang  (1990).  
The idea of the calculus is to suppose that all the data are sufficiently 
 regular and to expand all the
 variables with the Taylor formula.

\bigskip 
\vfill \eject 
\noindent  {\bf Proposition 3}. \quad Taylor expansion at zero  order.

 \noindent 
With the lattice Boltzmann defined previously, we have 

\smallskip \noindent   (5.1)  $ \qquad   \displaystyle
 f^j (x ,\, t) \,= \,  f^j_{\rm eq} (x ,\, t) \,+\, {\rm O} (\Delta t ) 
\,, \quad  0 \leq j  \leq J  \,, \, $ 

\smallskip \noindent   (5.2)  $ \qquad   \displaystyle
 f^j_* (x ,\, t) \,= \,  f^j_{\rm eq} (x ,\, t) \,+\, {\rm O} (\Delta t ) 
\,, \quad  0 \leq j  \leq J  \,, \, $ 

\smallskip \noindent  
with $\,  f^j_{\rm eq} \,$  defined from the conservative variables $ \, W \, $ 
according  to the  relation  (2.7). 

\smallskip \noindent  {\bf Proof of Proposition 3}.

 \noindent
The key point is to expand the relation (4.2) relative to the infinitesimal $ \, 
\Delta t .\,$ We have on one hand 

\smallskip \noindent   (5.3)  $ \qquad   \displaystyle
 f^j (x ,\, t+\Delta t ) \,= \,  f^j (x ,\, t)  \,+\, {\rm O} (\Delta t ) \, $ 

\smallskip \noindent 
and on the other hand 

\smallskip \noindent   (5.4)  $ \qquad   \displaystyle
 f^j_* (x -  v_j \,\Delta t ,\, t) \,= \,  f^j_* (x ,\, t)  \,+\, {\rm O} (\Delta t ) \, $ 

\smallskip \noindent 
Then  \quad  
$ \,  \displaystyle   m^{k}_* (x,\, t)  \, = \, \sum_{j=0}^{J} \, M^k_j \, f^j_*  (x,\, t) \,\,
 = \,\,  m^{k} (x,\, t)  \,+\, {\rm O} ( \Delta t )  \, $ \qquad  and 

\smallskip \noindent   (5.5)  $ \qquad   \displaystyle
  m^{k}_* (x,\, t)   \, - \,   m^{k} (x,\, t)   \, = \,  {\rm O} ( \Delta t )  \, . \, $  

\smallskip \noindent 
But, due  to (3.3), we have 

\smallskip \noindent   (5.6)  $ \qquad   \displaystyle
  m^{k}_* (x,\, t)   \, - \,   m^{k} (x,\, t)   \, = \, -s_k \, \big(  m^{k} (x,\, t) -
  m^{k}_{\rm eq} (x,\, t)  \big) \, .  \, $  

\smallskip \noindent  
From (5.5) and (5.6) we deduce, due to the fact that $ \, s_k \neq 0 \, $ when $  \, 
k \geq d+1 \, $: 

\smallskip \noindent   (5.7)  $ \qquad   \displaystyle
  m^{k} (x,\, t)   \, = \,   m^{k}_{\rm eq} (x,\, t)  +  {\rm O} ( \Delta t )  \, , \quad
  k \geq d+1 \,  . \, $ 

\smallskip \noindent
We insert (5.7) into (5.5) and we deduce 

\smallskip \noindent   (5.8)  $ \qquad   \displaystyle
  m^{k}_* (x,\, t)   \, = \,   m^{k}_{\rm eq} (x,\, t)  +  {\rm O} ( \Delta t )  \, , \quad
  k \geq d+1 \,  . \, $ 

\smallskip \noindent 
Taking into account the relations (2.14) and (3.1) on one hand and (2.12) and (3.9) on the
other hand, we deduce (5.1) and (5.2) from (5.7) and (5.8).   $\hfill \square$

\bigskip \bigskip    \noindent {\bf \large    6) \quad Taylor expansion at  first order  }    

\monitem
We expand now the relation (4.2) one step further with respect to 
 the time step $ \, \Delta t .\, $  We introduce the second order moment 

\smallskip \noindent   (6.1)  $ \qquad   \displaystyle
  F^{ \alpha \, \beta} \, \equiv \,   \sum_{j=0}^{J} \,
\,\,  v_j^{\alpha}  \,\,  v_j^{\beta} \,\,  f^{j}_{\rm eq}   \, , \quad
  1 \leq \alpha , \, \beta \leq d  \, . \, $ 

\smallskip \noindent 
We denote in the following $\, \partial_t \, $ instead of 
 $\, {{\partial}\over {\partial t}}  \, $ and  $\, \partial_{\beta} \, $
in  place of   $\, {{\partial}\over {\partial x_{\beta}}}  .\, $
Then we have the following result at the first order. 

\bigskip \noindent  {\bf Proposition 4}. \quad Euler equations of gas dynamics.

 \noindent 
With the lattice Boltzmann scheme previously defined, we have the conservation of mass and
momentum at the first order: 

\smallskip \noindent   (6.2)  $ \qquad   \displaystyle
  \partial_{t}  \rho   \,+\,   \sum_{\beta=1}^{d}   \partial_{\beta} \, q^{\beta}   \,\, = \,\,
{\rm O} ( \Delta t  )  \, $  

\smallskip \noindent   (6.3)  $ \qquad   \displaystyle
  \partial_{t}  q^{\alpha}    \,+\,  \sum_{\beta=1}^{d}   \partial_{\beta} \,  F^{ \alpha \, \beta}
\,\, = \,\,  {\rm O} ( \Delta t  )  \,. \,   $

\smallskip \noindent  {\bf Proof of Proposition 4}.

 \noindent
We expand  both sides of relation (4.2) up to first  order: 

\smallskip \noindent     $    \displaystyle
  f^{j}( x \,,\, t  \,+\, \Delta t)  \,  = \, 
f^{j}( x \,,\, t )  \, +\,   \Delta t \,\, \partial_t    f^{j}   \,  \,+\,
{\rm O} ( \Delta t^2  )  \, $

\smallskip \noindent     $    \displaystyle
   f^{j}_* ( x \,- \, v_j \,  \Delta t \,,\, t)  \,  = \, 
  f^{j}_* ( x \,,\, t)   \, - \,   \Delta t \,\, v_j^{\beta} \, \,  \partial_{\beta}
   f^{j}_*   \, \,+\, {\rm O} ( \Delta t^2  )  \, . \, $

\smallskip   \noindent 
We take the moment  of order $  \, k \, $ of this identity: 

\smallskip \noindent     $    \displaystyle
  m^k( x \,,\, t )  \,+\,  \Delta t \,\,  \partial_{t}   m^k   \,+\,
  {\rm O} ( \Delta t^2  )   \,= \,
  m^k_*( x \,,\, t )   \,-\,    \Delta t \,\, \sum_{j=0}^{J} \,
\, M^k_j  \,\,  v_j^{\beta} \,\, \partial_{\beta}   f^{j}_*   \,+\,
{\rm O} ( \Delta t^2  )  \, $

\smallskip   \noindent 
and   we use the previous Taylor expansions (5.1) (5.2)  at the order zero: 

\smallskip \noindent   (6.4)  $ \quad   \displaystyle
  m^k( x \,,\, t )  \,+\,   \Delta t \,\,  \partial_{t}   m^k_{\rm eq}   \, = \,
    m^k_*( x \,,\, t )   \,-\,    \Delta t \, \,   \sum_{j=0}^{J} \,
\,\, M^k_j  \,\,  v_j^{\beta} \,\, \partial_{\beta}   f^{j}_{\rm eq}   \,+\,
{\rm O} ( \Delta t^2  )  \, . \, $

\smallskip \noindent 
We take $ \, k=0 \, $ inside the  relation (6.4). We get (6.2) since 
$ \,   m^0 ( x \,,\, t ) \equiv      m^0_* ( x \,,\, t )  \equiv \rho  ( x \,,\, t ) .\, $ 
Considering now the particular case $ \, k=\alpha \, $ with 
$  \, 1 \leq \alpha \leq d ,\, $ we have also  $ \,   m^\alpha( x \,,\, t ) \equiv 
    m^\alpha_*( x \,,\, t )   \equiv q^\alpha  ( x \,,\, t ) \, $ 
and the relation (6.3) is a direct consequence of the
    definition (6.1) and the property (2.11).      $\hfill \square$ 

\bigskip \noindent  {\bf Proposition 5}. \quad Technical lemma.

 \noindent 
We introduce the ``conservation defect'' $ \, \theta ^k \, $ according to the relation 

\smallskip \noindent   (6.5)  $ \quad   \displaystyle 
\theta^k  ( x \,,\, t ) \,=\,  \partial_{t}   m^k_{\rm eq} \, +\, \sum_{j=0}^{J}
  \, M_j^k \,\, v_j^{\beta} \,\, \partial_{\beta}    f^{j}_{\rm eq}
\, \equiv \,  \sum_{j=0}^{J}
\, M_j^k \,\, ( \, \partial_{t}   f^{j}_{\rm eq} \,+\,  v_j^{\beta} \,\,
\partial_{\beta}    f^{j}_{\rm eq}   \, ) \, . \, $

\smallskip \noindent 
Then we have the following properties: 

\smallskip \noindent   (6.6)  $ \quad   \displaystyle 
 m^{k}( x \,,\, t )   \, = \, m^k_{\rm eq}( x \,,\, t )  \,-\, {{\Delta t }\over {s_k}} \,\,
\theta^k   \,  \,+\, {\rm O} ( \Delta t^2 ) \,, \quad k \geq d+1  \,,  \, $

\smallskip \noindent   (6.7)  $ \quad   \displaystyle  
m^{k}_* ( x \,,\, t )  \, = \, m^k_{\rm eq}( x \,,\, t )  \,-\, \Big( {{1}\over{s_k}} \,-\, 1
\Big) \,\, \Delta t  \,\, \theta^k   \,  \,+\, {\rm O} ( \Delta t^2 )  
 \,, \quad k \geq d+1  \,,  \, $

\smallskip \noindent   (6.8)  $ \quad   \displaystyle  
\partial_{\beta}   f^{j}_*   \, = \,   \partial_{\beta}   f^{j}_{\rm eq}
 \,-\, \Delta t  \,\, \sum_{k = d+1}^{J} \,  \Big( {{1}\over{s_k}} \,-\, 1 \Big) \,\, (M^{-1})_k^j \,\,
 \partial_{\beta}  \theta^k   \,  \,+\, {\rm O} ( \Delta t^2 )  \, .\, $

\smallskip \noindent  {\bf Proof of Proposition 5}.

 \noindent
We start from the relation (6.4) and we have observed at the previous proposition that 

\smallskip \noindent   (6.9)  $ \quad   \displaystyle  
  \theta^i  \,=\, {\rm O} ( \Delta t ) \, , \quad 0 \leq i \leq d \, . \, $  

\smallskip \noindent 
We remark also that from the relation (5.6), we have 

\smallskip \noindent   $    \displaystyle  
m^k ( x \,,\, t ) - m^k_{\rm eq} ( x \,,\, t ) \, = \, {{1 }\over{ s_k }} \, \big( 
 m^{k} ( x \,,\, t ) \, - \,  m^{k}_*  ( x \,,\, t ) \big)   \quad   {\rm if }   \,\, k \geq d+1 .\,  $

\smallskip \noindent 
Then the relation (6.6) is a direct consequence of (6.4) and the definition (6.5). In
consequence, the relation (6.7) follows from (6.6) and (6.4). Due to (6.7),   (6.9) and (3.9), we
have 

\smallskip \noindent   (6.10)  $ \quad   \displaystyle  
  f^{j}_* ( x \,,\, t )   \, = \,    f^{j}_{\rm eq} ( x \,,\, t )
 \,-\, \Delta t  \,\, \sum_{k \geq d+1} \,  \Big( {{1}\over{s_k}} \,-\, 1 \Big) \,\, (M^{-1})_k^j \,\,
  \theta^k   \,  \,+\, {\rm O} ( \Delta t^2 )  \, $

\smallskip \noindent
and the relation (6.8) follows from derivating (6.10) in the  direction 
 $ \,x_{ \beta} .  \, $   $\hfill \square$ 

\bigskip \bigskip    \noindent {\bf \large    7) \quad Equivalent  equation at  second  order }    

\monitem
We  introduce the tensor $\,  \Lambda ^{\alpha \, \beta}_k \, $ according to 

\smallskip \noindent   (7.1)  $ \quad   \displaystyle  
\Lambda ^{\alpha \, \beta}_k \,\, \equiv \,\,
  \sum_{j=0}^{J} \,  v_j^{ \alpha} \,\, v_j^{\beta} \,\, (M^{-1})_k^j  \, , \quad
  1 \leq \alpha , \, \beta \leq d   \, , \quad  0 \leq k  \leq J  \,. \, $ 

\smallskip \noindent  
We can now establish the major result of our contribution. 

\bigskip \noindent  {\bf Proposition 6}. \quad Navier-Stokes equations of gas dynamics.

 \noindent 
With the lattice Boltzmann method defined in previous sections 
and the conservation defect $ \, \theta ^k \, $ defined in (6.5), 
we have the following
expansions up to second order accuracy: 

\smallskip \noindent   (7.2)  $ \qquad   \displaystyle
  \partial_{t}  \rho   \,+\,   \sum_{\beta=1}^{d}   \partial_{\beta} \, q^{\beta}   \,\, = \,\,
{\rm O} ( \Delta t^2  )  \, $  

\smallskip \noindent   (7.3)  $ \qquad   \displaystyle
\partial_{t}  q^{ \alpha}   \,+\,   \sum_{\beta=1}^{d} \,   \partial_{\beta}  \,
\bigg( F ^{\alpha \, \beta} -  \Delta t   \sum_{k \geq d+1} \,
 \Big( {{1}\over{s_k}} - {1\over2}  \Big)  \,  \Lambda ^{\alpha \, \beta}_k  \,\,
 \,   \theta^k  \bigg)  \, = \, {\rm O} ( \Delta t ^2 ) \, .  \, $ 

\monitem
A consequence of relation (7.3) is the fact that a lattice Boltzmann scheme approximates 
at second order of accuracy a Navier-Stokes type equation with viscosities $\, \mu_k\, $ of
the form

\smallskip \noindent   (7.4)  $ \qquad   \displaystyle
\mu_k   \, = \, \Delta t \,  \Big( {{1}\over{s_k}} - {1\over2}  \Big) \, .\, $ 

\smallskip \noindent 
We refer for the details to D. D'Humi\`eres (1992), Lallemand and Luo (2000) or
to our recent survey (2007). The relations (7.4) are  known as the ``D'Humi\`eres
relations''.  We observe that in practice, the scalar  $\, \mu_k\, $ is imposed by the
physics and by the parameter   $\,  \Delta t \, $ is constrained by the space discretization
  $\,  \Delta x \, $
and the relation (2.1). Then the parameter  $\, s_k\, $ must  be chosen in order to satisfy
the D'Humi\`eres relations (7.4).

\smallskip \noindent  {\bf Proof of Proposition 6}.

 \noindent
We start again from the identity (4.2). We expand both terms up to second order accuracy:

\smallskip \noindent   $    \displaystyle 
 f^{j}( x \,,\, t  \,+\, \Delta t)  \,=\,  
f^{j}( x \,,\, t )  \, \,+\, \,  \Delta t \,\, \partial_t    f^{j}   \,  \,+\,
{{1}\over{2}} \,  \Delta t ^2  \,\,  \partial_{tt}^2    f^{j}   \,  \,+\,
{\rm O} ( \Delta t^3  )  \, $

\smallskip \noindent   $    \displaystyle 
  f^{j}_* ( x \,- \, v_j \,  \Delta t \,,\, t)  \, = \,  
  f^{j}_* ( x \,,\, t)    \,-\,  \,   \Delta t \,\, v_j^{\beta} \,\,
 \partial_{\beta}  f^{j}_*   \, \,+\,
{{1}\over{2}} \,  \Delta t ^2  \,\,   v_j^{\beta} \,  v_j^{\gamma} \,\,  
\partial_{\beta  \gamma }^2  f^{j}_*     \,+\,  {\rm O} ( \Delta t^3  )  \, . \, $

\smallskip \noindent  
We take the moment of order $  \, i \,\, (0 \leq i \leq  d) \,  \, $ of this
identity. We obtain:

\setbox21=\hbox{$\displaystyle  
m^i ( x \,,\, t)    \,+\,  \Delta t \,    \,  \partial_{t}   m^i   \,+\,
{{1}\over{2}} \,  \Delta t ^2  \,\,  \partial_{tt}^2     m^i  \,  \,+\,
{\rm O} ( \Delta t^3  )   \,\, = \,\,  m^i_* ( x \,,\, t)  \,+\,    $}
\setbox22=\hbox{$\displaystyle \quad 
 \,-\,    \Delta t \,\, \sum_{j=0}^{J}
\,\, M^i_j  \,\,  v_j^{\beta} \,\, \partial_{\beta}   f^{j}_*  
\,+\, {{1}\over{2}} \,  \Delta t ^2  \, \sum_{j=0}^{J} \, M^i_j  \,
v_j^{\beta} \,  v_j^{\gamma} \,\,   \partial_{\beta  \gamma }^2  f^{j}_*  
\,+\, {\rm O} ( \Delta t^3  )  \, . \,   $}
\setbox30= \vbox {\halign{#\cr \box21 \cr   \box22 \cr    }}
\setbox31= \hbox{ $\vcenter {\box30} $}
\setbox44=\hbox{\noindent  (7.5) $\displaystyle \, \,   \left\{ \box31 \right. $}  
\smallskip \noindent $ \box44 $

\smallskip \noindent  
We use the microscopic conservation $ \,\,   m^i_* ( x \,,\, t)   \equiv  m^i ( x \,,\, t)
\,\,$ in (7.5) and   the previous Taylor expansion  at  order one, 
in particular the relation (6.8).  We divide by $\, \Delta t \, $ and we deduce: 

\smallskip \newpage \noindent   $    \displaystyle 
 \partial_{t}   m^i   \,+\,
{{1}\over{2}} \,  \Delta t  \,\,  \partial_{tt}^2     m^i  \, = \,
  -   \sum_{j=0}^{J}  \, M^i_j  \,\,  v_j^{\beta} \,\,
 \partial_{\beta}   f^{j}_{\rm eq} \,\,+\,  $

\smallskip  \noindent   $    \displaystyle \qquad  \qquad \qquad \qquad \qquad
 \,  +  \,  \Delta t  \,  \,  \sum_{j=0}^{J}  \, \sum_{k \geq d+1}
\! M^i_j  \,  v_j^{\beta} \,  \Big( {{1}\over{s_k}} - 1 \Big) \, (M^{-1})_k^j \,\,
 \partial_{\beta}   \, \theta^k   \,  \,+\, $ 

\smallskip  \noindent   $    \displaystyle \qquad   \qquad \qquad \qquad \qquad \qquad
 \,  +  \, {{1}\over{2}} \,  \Delta t  \,  \,  \sum_{j=0}^{J}  \,\, M^i_j  \,\,  
v_j^{\beta} \,  v_j^{\gamma} \,\,   \partial_{\beta  \gamma }^2  f^{j}_{\rm eq}
\,+\, {\rm O} ( \Delta t^2  )  \, . \, $

\noindent 
Then

\setbox21=\hbox{$\displaystyle  
 \partial_{t}   m^i   \,+\, \sum_{\beta=1}^{d}  \, \sum_{j=0}^{J} 
 \, M^i_j  \,\,  v_j^{\beta} \,\,
 \partial_{\beta}   f^{j}_{\rm eq}   \,\, = \,\,     $}
\setbox22=\hbox{$\displaystyle \qquad   \,\, = \,\,  
 \Delta t \,  \sum_{\beta=1}^{d}  \,   \sum_{j=0}^{J}    \, \sum_{k \geq d+1}
\, M^i_j  \,\,  v_j^{\beta} \,  \Big( {{1}\over{s_k}} - 1 \Big) \, (M^{-1})_k^j \,\,
 \partial_{\beta}  \theta^k   \,+\,  \,   $}
\setbox23=\hbox{$\displaystyle  \qquad   \qquad  \,\,
 \,+\, \,   {{\Delta t }\over{2}} \, \Big( -  \partial_{tt}^2     m^i  \,  \,+\,
 \,\  \sum_{\beta=1}^{d}  \,  \sum_{j=0}^{J}    \,\, M^i_j  \,\,    v_j^{\beta} \,  
v_j^{\gamma} \,\,  
 \partial_{\beta  \gamma }^2  f^{j}_{\rm eq} \, \Big)
\,+\, {\rm O} ( \Delta t^2  )  \, . \, $}
\setbox30= \vbox {\halign{#\cr \box21 \cr   \box22 \cr \box23 \cr    }}
\setbox31= \hbox{ $\vcenter {\box30} $}
\setbox44=\hbox{\noindent  (7.6) $\displaystyle \quad    \left\{ \box31 \right. $}  
\smallskip  \noindent $ \box44 $

\monitem
We set    $ \,  i = 0 \,$  in  the relation (7.6) and we  look for  the
conservation of mass. Due to the property $ \, M^0_j  \equiv  1 , \,    $  the  sum over
$ \,  j \, $  in the second line of (7.6) is null since 
$ \, \smash {   \sum_{j=0}^{J}  }  v_j^{\beta} \,  (M^{-1})_k^j \, $
is equal to zero. 
We have also the following algebraic calculus: 

\smallskip \noindent   $    \displaystyle 
\partial_{tt}^2    m^0 \, =\,
\partial_{tt}^2    \rho \, =\,  -   \sum_{\beta=1}^{d}  \,  \partial_{t \beta }^2\,
q^{\beta}  \,+\, {\rm O} ( \Delta t  ) 
\, =\,   -  \sum_{\beta=1}^{d}  \,  \partial_{\beta} \, \partial_{t}  
\, q^{\beta}  \,+\, {\rm O} ( \Delta t  )  \, =\,   $ 

\smallskip \noindent   $    \displaystyle \qquad \qquad \qquad
  \, =\, 
 \sum_{\beta=1}^{d}  \,  \sum_{\gamma=1}^{d}  \,
 \partial_{\beta \gamma}^2  \,  F^{  \beta  \,\gamma  }   \,+\, {\rm O} ( \Delta t  )  
 \, =\,   \sum_{\beta=1}^{d}  \,  \sum_{\gamma=1}^{d}  \, \, \sum_{j=0}^{J}   
 \,\,    v_j^{\beta} \,  v_j^{\gamma} \,\,  
 \partial_{\beta  \gamma }^2  f^{j}_{\rm eq}  \,+\, {\rm O} ( \Delta t  )  \,  $

\smallskip \noindent  
and the third line of  (7.6) is null up to second order accuracy.
  Thus  the conservation of mass (7.2) up to second order  accuracy is established. 

\monitem
We set    $ \,  i =  \alpha \,$ with   $ \, 1 \leq \alpha \leq d \, $ and we look for
 the   conservation of momentum. In this particular case, the relation (7.6) takes the
form:

\setbox21=\hbox{$\displaystyle  
 \partial_{t}  q^{ \alpha}   \,+\,  \sum_{\beta=1}^{d}  \, \sum_{j=0}^{J} \,  
    v_j^{ \alpha}  \,\,  v_j^{\beta} \,\,
 \partial_{\beta}   f^{j}_{\rm eq}   \,\, = \,\,     $}
\setbox22=\hbox{$\displaystyle \qquad  \quad   \,\, = \,\,  
  \Delta t    \sum_{k \geq d+1}
\,  \Big( {{1}\over{s_k}} - 1 \Big) \, \sum_{\beta=1}^{d}   \, \Big[ \,
\sum_{j=0}^{J}  \,    v_j^{ \alpha}  \,\, v_j^{\beta} \,\, (M^{-1})_k^j \,\Big] \,
 \partial_{\beta}  \theta^k   \,+\,   $}
\setbox23=\hbox{$\displaystyle  \qquad  \qquad 
\,+\, {{\Delta t }\over{2}} \, \Big( -  \partial_{tt}^2   q^{ \alpha} \,  \,+\,
 \,\, \sum_{\beta=1}^{d}  \,     \sum_{j=0}^{J} \,     v_j^{ \alpha}  \,
   v_j^{\beta} \,  v_j^{\gamma} \,\,  
 \partial_{\beta  \gamma }^2  f^{j}_{\rm eq} \, \Big)
\,+\, {\rm O} ( \Delta t^2  )  \,. \,   $}
\setbox30= \vbox {\halign{#\cr \box21 \cr   \box22 \cr \box23 \cr    }}
\setbox31= \hbox{ $\vcenter {\box30} $}
\setbox44=\hbox{\noindent  (7.7) $\displaystyle \quad    \left\{ \box31 \right. $}  
\smallskip \noindent $ \box44 $

\smallskip \noindent
We have now to play with some    algebra: 

\smallskip
\newpage \noindent   $    \displaystyle 
 -  \partial_{tt}^2   q^{ \alpha} \,  \,+\,
 \,  \sum_{\beta=1}^{d}  \,   \sum_{j=0}^{J} \,  \,   v_j^{ \alpha}  \,
    v_j^{\beta} \,  v_j^{\gamma} \,\,  
 \partial_{\beta  \gamma }^2  f^{j}_{\rm eq}  \,\,=\,\, $ 

\smallskip \noindent   $    \displaystyle \qquad \qquad \,\,=\,\, 
 \sum_{\beta=1}^{d}  \,  \Big( \,   \partial_{t}   \partial_{\beta} F^{ \alpha \, \beta} +
 \,  \sum_{j=0}^{J} \,  \,   v_j^{ \alpha}  \,\,    v_j^{\beta} \,  v_j^{\gamma} \,\,  
 \partial_{\beta  \gamma }^2  f^{j}_{\rm eq}   \,  \Big)  \,+\, {\rm O} ( \Delta t  )    \,$

\smallskip \noindent   $    \displaystyle \qquad \qquad \,\,=\,\, 
 \sum_{\beta=1}^{d}  \,   \partial_{\beta} \,  \Big( \,
 \,  \sum_{j=0}^{J} \,  \,   v_j^{ \alpha}  \,\,    v_j^{\beta} \, \big(
 \partial_{t}  f^{j}_{\rm eq} \,+\,  v_j^{\gamma} \, \partial_{\gamma }  f^{j}_{\rm eq}
 \big) \, \Big)  \,+\, {\rm O} ( \Delta t  )  \, $

\smallskip \noindent   $    \displaystyle \qquad \qquad \,\,=\,\, 
 \sum_{\beta=1}^{d}  \,   \partial_{\beta} \,  \Big( \,
 \,  \sum_{j=0}^{J} \,  \,   v_j^{ \alpha}  \,\,    v_j^{\beta} \,
\,\,  \sum_{k=0}^{J} \, (M^{-1})^j_k \, \theta^k  \, \Big)  \,+\, {\rm O} ( \Delta t  ) \, $

\smallskip \noindent   $    \displaystyle \qquad \qquad \,\,=\,\, 
 \sum_{\beta=1}^{d}  \,   \partial_{\beta} \,  \Big( \,
 \,\,\sum_{k \geq d+1}   \,   \Big[  \sum_{j=0}^{J} \, 
  v_j^{ \alpha}  \,\, v_j^{\beta} \,\, (M^{-1})_k^j \,\Big] \,
\, \theta^k \,   \, \Big) \,  \,+\, {\rm O} ( \Delta t  )   \, $

\smallskip \noindent   $    \displaystyle \qquad \qquad \,\,=\,\, 
 \sum_{\beta=1}^{d}  \,   \partial_{\beta} \,  \Big( \,
 \,\,\sum_{k \geq d+1}   \,  \Lambda ^{\alpha \, \beta}_k \,
\, \theta^k \,   \, \Big) \,  \,+\, {\rm O} ( \Delta t  )  \, $

\smallskip \noindent 
due to the definition (7.1). We deduce from (6.1), (7.7) and the above calculus: 

\smallskip \noindent   $    \displaystyle 
\partial_{t}  q^{ \alpha}   \,+\,  \sum_{\beta=1}^{d}  \,   
\partial_{\beta} F ^{\alpha \, \beta}  \,\,=\,\,   \Delta t    \sum_{k \geq d+1}
\,  \Big( {{1}\over{s_k}} - 1 \Big) \, \sum_{\beta=1}^{d}  \, \Lambda ^{\alpha \, \beta}_k \,
  \,  \, \partial_{\beta}   \, \theta^k   \,+\,   $

\smallskip \noindent   $    \displaystyle \qquad \qquad  \qquad   \qquad   \qquad 
\,\,+\,\,   {{\Delta t }\over{2}} \,  \, 
 \sum_{\beta=1}^{d}  \,   \partial_{\beta} \,  \Big( \,
 \,\,\sum_{k \geq d+1}   \,  \Lambda ^{\alpha \, \beta}_k \,
\, \theta^k \,   \, \Big) \,  \,+\, {\rm O} ( \Delta t ^2  )  \, $

\smallskip \noindent   $    \displaystyle \qquad \qquad  \qquad   \qquad  \,\,=\,\,  
\Delta t \,   \sum_{\beta=1}^{d}  \,  \sum_{k \geq d+1} \,
 \Big( {{1}\over{s_k}} - {1\over2}  \Big)  \,  \Lambda ^{\alpha \, \beta}_k  \,\,
 \partial_{\beta}  \, \theta^k  \, + \, {\rm O} ( \Delta t ^2 ) \, .  \, $ 

\smallskip \noindent 
and the relation (7.3) is established.    $\hfill \square$

\bigskip \bigskip    \noindent {\bf \large    8) \quad Equivalent  equation at  second  order }    

\monitem
The previous propositions establish that the equivalent partial differential 
 equations of a Boltzmann scheme are 
given up to second order accuracy by the same result as the formal Chapman-Enskog
expansion. We find Euler type equation at the first order (Proposition~4) and
Navier-Stokes type equation at the second order (Proposition 6). 
Note that with the above framework  no  {\it a priori} formal two-time multiple 
 scaling  is necessary to establish the  Navier-Stokes
 equations from a lattice Boltzmann scheme, as  done previously in the contribution
of D'Humi\`eres.    We remark also that a so-called diffusive scaling like 
$\, \smash { {{\Delta t}\over{\Delta x^2}} }  =  $ 
constant, instead of our condition (2.1)  $\, {{\Delta t}\over{\Delta x}} = $ 
 constant,  leads   to the incompressible Navier-Stokes equations, as proposed  by Junk,
 Klar and Luo (2005). 
In both cases, we have just to use  the
 Taylor formula for  a single infinitesimal parameter.

\bigskip \bigskip   \newpage \noindent {\bf \large    9) \quad Acknowledgments }    

\noindent
The author thanks Li-Shi Luo to his kind invitation to present the scientific work of 
Orsay's team at ICMMES Conference 
in July 2005.   The author thanks also with a great emphasis   
Pierre Lallemand for very  helpfull  discussions all along the elaboration 
of this contribution. 
Last but not least, the referees transmitted to the author very good remarks that have
been incorporated inside  the present edition of the article.

\bigskip \bigskip    \noindent {\bf \large    10) \quad References }    



\smallskip 
 \hangindent=7mm \hangafter=1 \noindent 
S.C. Chang. ``A critical analysis of the modified equation technique of Warming and
Hyett'',  {\it Journal of Computational Physics}, vol. 86, p.~107-126, 1990. 

 \hangindent=7mm \hangafter=1 \noindent   
D. D'Humi\`eres. ``Generalized Lattice-Boltzmann Equations'', 
in {\it Rarefied Gas Dynamics: Theory
and Simulations}, vol. 159 of {\it AIAA Progress in
Astronautics and Astronautics}, p. 450-458, 1992.

 \hangindent=7mm \hangafter=1 \noindent   
F. Dubois. ``Une introduction au sch\'ema de Boltzmann sur r\'eseau'', 
{\it ESAIM Proceedings}, vol.~18, pages 181-215, July 2007.

 \hangindent=7mm \hangafter=1 \noindent   
D. Griffiths, J. Sanz-Serna. ``On the scope of the method of modified equations'', 
{\it SIAM Journal on  Scientific and Statistical  Computing}, vol.~7, p.~994-1008, 1986.

 \hangindent=7mm \hangafter=1 \noindent   
 M. Junk, A. Klar, and L.-S. Luo. ``Asymptotic analysis of the lattice 
Boltzmann equation'', {\it 
Journal of Computational Physics}, vol.  210, p.~676-704, December 2005.

\smallskip \hangindent=7mm \hangafter=1 \noindent  
P. Lallemand, L.-S. Luo. ``Theory of the lattice Boltzmann method: 
Dispersion, dissipation, isotropy, Galilean invariance, and stability'',
{\it Physical Review E}, vol.  61, p. 6546-6562, June 2000. 
 
\smallskip \hangindent=7mm \hangafter=1 \noindent  
A. Lerat, R. Peyret, ``Noncentered Schemes and Shock Propagation
  Problems'',  {\it  Computers and Fluids}, vol. {\bf 2}, p. 35-52, 1974. 

\smallskip \hangindent=7mm \hangafter=1 \noindent   
G. Strang. {\it An introduction to applied mathematics}, Wellesley-Cambridge press, 
 Wellesley,  1986. 

\smallskip \hangindent=7mm \hangafter=1 \noindent   
R.F. Warming,  B.J. Hyett, ``The modified equation approach 
to the stability and accuracy analysis of finite difference methods'', 
{\it Journal of Computational Physics}, vol. {\bf 14}, p. 159-179, 1974.


\end{document}